\documentclass{amsart}

\usepackage{amssymb}

\newtheorem{theorem}{Theorem}[section]

\theoremstyle{definition}
\newtheorem{definition}[theorem]{Definition}

\theoremstyle{corollary}
\newtheorem{corollary}[theorem]{Corollary}
\theoremstyle{proposition}
\newtheorem{proposition}[theorem]{Proposition}

\theoremstyle{remark}
\newtheorem{remark}[theorem]{Remark}

\numberwithin{equation}{section}

\def\Dc{{\mathcal D}}

\def\e{{\rm e}}

\def\eh{{{\bf{\rm{e}}}}}

\def\Dc{{\mathcal D}}

\def\eh{{{\bf{\rm{e}}}}}

\def\Dc{{\mathcal D}}

\def\eh{{{\bf{\rm{e}}}}}

\def\Dc{{\mathcal D}}

\begin{document}

\title{The ADM mass of asymptotically flat hypersurfaces}


\author{Levi Lopes de Lima}
\address{Federal University of Cear\'a,
Department of Mathematics, Campus do Pici, Av. Humberto Monte, s/n, 60455-760,
Fortaleza/CE, Brazil.}
\curraddr{}
\email{levi@mat.ufc.br}
\thanks{Partially supported by CNPq/BR and FUNCAP/CE}

\author{Frederico Gir\~ao}
\address{Federal University of Cear\'a,
Department of Mathematics, Campus do Pici, Av. Humberto Monte, s/n, 60455-760,
Fortaleza/CE, Brazil.}
\curraddr{}
\email{fred@mat.ufc.br}
\thanks{}


\date{}

\dedicatory{}

\begin{abstract}
We provide integral formulae for the ADM mass of asymptotically flat hypersurfaces in Riemannian manifolds with a certain warped product structure in a neighborhood of \lq spatial\rq\, infinity, thus extending Lam's recent results   on Euclidean graphs \cite{L1} \cite{L2} to this broader context. As applications we exhibit, in {\em any} dimension, new examples of manifolds for which versions of the Positive Mass and Riemannian Penrose inequalities hold and discuss a notion of quasi-local mass in this setting. The proof explores a novel connection between the co-vector  defining the ADM mass of  a hypersurface as above and the Newton tensor associated to its shape operator, which takes place in the presence of an ambient Killing field.
\end{abstract}

\maketitle


\section{Introduction and statement of results}
\label{intstate}

An {\em asymptotically flat end} is a  Riemannian manifold $(E, h)$ of dimension $n\geq 3$ for which  there exists a diffeomorphism $\Psi:\mathbb R^n-\overline B_1(0)\to E$ introducing coordinates in $E$, say $x=(x_1,\cdots,x_n)$, such that  the following decay conditions hold as $|x|\to +\infty$ for some $\tau>(n-2)/2$:
\begin{equation}\label{decay1}
h_{ij}=\delta_{ij}+O(|x|^{-\tau}),
\quad h_{ij,k}=O(|x|^{-\tau-1}).
\end{equation}
Here, the $h_{ij}$'s are the coefficients of $h$ with respect to $x$, $h_{ij,k}=\partial h_{ij}/\partial x_k$ and $|\,\,|$ is the standard Euclidean norm. We also assume that the scalar curvature $R_h$ of  $h$ is integrable.
Under these conditions the {\em ADM mass} of $(E,h)$ is defined by
\begin{equation}\label{mass}
m_h=\lim_{r\to +\infty}c_n\int_{\Sigma_r}(h_{i,j}^j-h_{j,i}^j)\nu^id\Sigma_r,\quad c_n=\frac{1}{2(n-1)\omega_{n-1}},
\end{equation}
where $\nu=(\nu^1,\cdots, \nu^n)$ is the outward unit normal to a large coordinate sphere $\Sigma_r$ of radius $r$ and $\omega_{n-1}$ is the volume of the unit sphere of dimension $n-1$. Unless otherwise stated, all manifolds considered here are connected and smooth. Also, throughout the paper we are summing over repeated indices.

Despite being phrased in local coordinates at infinity, it is known that the mass is well-defined and finite under the above decay conditions \cite{Ba1} \cite{Ch} \cite{M}. This invariant, which had its origins in the  Hamiltonian formulation of General Relativity \cite{ADM}, proved itself extremely useful in Geometric Analysis. It plays, for instance, a crucial role in questions of existence and compactness of solutions for the Yamabe problem \cite{S} \cite{LP} \cite{LZ} \cite{KMS} \cite{BM}.

The {\em rationale} behind the concepts above is of course the Positive Mass Conjecture.
Recall that a {complete} Riemannian manifold $(M,g)$ is {\em asymptotically flat} if there exists a compact subset $K\subset M$ such that $E_M=M-K$, the end of $M$, is asymptotically flat in the restricted metric. By definition, the mass of $(M,g)$ is the mass $m_g$ of $(E_M,g)$.  The conjecture then
says that $R_g\geq 0$ implies $m_g\geq 0$, with equality holding if and only if $(M,g)=(\mathbb R^n,g_0)$, where $g_0$ is the standard flat metric. This assertion has been first proved by Schoen and Yau if $n\leq 7$ \cite{SY} and subsequently a proof for spin manifolds, which works in any dimension, was provided by Witten \cite{W} \cite{Ba1} \cite{LP}.
Recently, Lam \cite{L1} \cite{L2} gave an elementary proof of the mass inequality in the case that $(M,g)$ can be isometrically embedded as a complete graph in $(\mathbb R^{n+1},g_0)$; see also \cite{HW}, where a related argument is presented for the case of embedded hypersurfaces. Our intention here is to revisit Lam's argument and discuss a few generalizations.

In order to put our results in their proper perspective, one should discuss the contributions above in more detail. First, Schoen  and Yau established their result by assuming $R_g\geq 0$ and $m_g<0$ and then using minimal surfaces techniques to reach a contradiction, a method that does not  directly relate the energy density of the gravitational system modeled by $(M,g)$,  represented by the scalar curvature, to its mass as measured at spatial infinity by (\ref{mass}). On the other hand, Witten is able to express the mass as an integral of a quantity which is manifestly nonnegative if $R_g\geq 0$, but his formula is somewhat mysterious since it involves the choice of a suitable spinor  on $M$. Lam's approach in its turn, even though restricted to Euclidean graphs,  effectively  relates in a simple way the data of the problem through the remarkable formula
\begin{equation}\label{lam}
m_g=c_n\int_M\frac{R_g}{\sqrt{1+|df|^2}}dM,
\end{equation}
where $f$ is the function defining the graph.
Lam's reasoning is based on a divergence type formula for the scalar curvature of graphs in nonparametric coordinates reminiscent of an old result due to Reilly \cite{R}, so that (\ref{lam}) follows by a simple integration by parts argument. Here, we essentially follow the same strategy but our starting point is instead a closely related flux-type formula (see (\ref{basic}) below) for hypersurfaces (not necessarily graphs) in the presence of a Killing field, a result that appears in various guises in the literature; see, for instance, \cite{ABC} \cite{AdLM} \cite{AM} \cite{ARS} \cite{dL} \cite{Ro}. As confirmed by the arguments put forward below, this formula happens to provide a notable connection between the intrinsic geometry at infinity of an asymptotically flat hypersurface, as captured by the ADM mass, and its extrinsic geometry as encoded in the Newton tensor (see (\ref{zero}) below) of the shape operator; see Remark \ref{novel} below. Here we illustrate this principle by showing how Lam's formula can be extended to asymptotically flat hypersurfaces in certain Riemaniann manifolds with a warped product structure in a neighborhood of \lq spatial\rq\, infinity. One should remark however that this method does not seem to be well suited to address the important issue of rigidity, which is so useful in applications.

Our approach also extends the Penrose-like formula in \cite{L1} to this broader context. Recall that the (Riemannian) Penrose inequality is a conjectured sharpening of the positive mass inequality when the asymptotically flat manifold is allowed to have a compact \lq inner\rq\, boundary $\Gamma$. Thus, if $\Gamma$ is a (possibly disconnected) outermost minimal hypersurface of area $A$, the inequality says that
\begin{equation}\label{pen}
m_g\geq\frac{1}{2}\left(\frac{A}{\omega_{n-1}}\right)^{\frac{n-2}{n-1}},
\end{equation}
whenever $R_g\geq 0$, with equality holding if and only if $(M,g)$ is the Riemannian Schwarzschild solution.
Physically, $\Gamma$ is thought of as being the horizon of a collection of black holes whose overall contribution to the mass is through its total surface area as indicated  above.

If $n=3$ the conjecture has been verified for $\Gamma$ connected by Huisken and Ilmanen \cite{HI} and in general by Bray \cite{Br}. More recently, Bray and Lee \cite{BL} established the conjecture for $n\leq 7$ with the extra requirement that $M$ be spin for the rigidity statement. Even though many partial results have been obtained \cite{H1} \cite{Sc} \cite {FS} \cite{J}, the conjecture remains wide open in higher dimensions except for the case of Euclidean graphs recently investigated by Lam \cite{L1}. Thus, if  $(M,g)\subset (\mathbb R^{n+1},g_0)$ is an asymptotically flat graph with an inner boundary $\Gamma$ whose connected components lie on (possibly distinct)  horizontal hyperplanes and if we further assume that $M$ is orthogonal to the hyperplanes along $\Gamma$, it is proved in \cite{L1} that
\begin{equation}\label{penrmassger}
m_g=c_n\int_{\Gamma}s_1 d\Gamma+c_n\int_M\frac{R_g}{\sqrt{1+|df|^2}}\,dM,
\end{equation}
where $s_1$ is the mean curvature of  $\Gamma$, viewed as a planar hypersurface . In particular,  if $R_g\geq 0$ and each connected component of $\Gamma$ is convex, (\ref{pen}) follows by from the so-called Alexsandrov-Fenchel inequality as explained in \cite{L1}. It turns out that this circle of ideas can be considerably generalized, as we now pass to describe. As a  consequence we will exhibit, in any dimension, new examples of manifolds for which versions of the Positive Mass and  Penrose inequality hold.
Also, as  a by-product of our computation, we discuss an extension to this more general setting of a notion of quasi-local mass first considered in \cite{L2}.

We now explain our setup. We consider an  $n$-dimensional asymptotically flat end $(E,h)$. As usual, we fix  asymptotically flat coordinates $x=(x_1,\cdots,x_n)$ on $E$  by means of a diffeomorphism $\Psi:\mathbb R^n-\overline B_1(0)\to E$ and set $\eh_i=\partial/\partial x_i$, so that the coefficients of $h$ in these coordinates are $h_{ij}=\langle \eh_i,\eh_j\rangle$, where $\langle\,,\,\rangle$ is the inner product associated to $h$. These coefficients satisfy the corresponding decay conditions at infinity as in (\ref{decay1}).
We also fix a {\em positive} smooth function $\phi:E\to \mathbb R$
satisfying
\begin{equation}\label{um}
\lim_{r\to +\infty}\|\phi\circ\Psi-1\|_{C^1(\mathbb R^n-B_r(0))}=0,
\end{equation}
and consider the warped product $(\overline{E},\overline h)$, with ${{\overline E}}=E\times {{I}}$ and $\overline h=h+\phi^2dt^2$, where $t$ is the standard linear coordinate in $I\subset \mathbb R$. Notice that each $t\in I$ defines a horizontal slice $E_t=E\times \{t\}\hookrightarrow \overline E$ which is totally geodesic, so that $E_t=E$ isometrically. This follows easily from the fact that $X=\partial/\partial t$, the vertical coordinate field, is Killing. Notice, moreover, that from $\phi=|X|_{\overline h}$ we find that
\begin{equation}\label{framevert}
\eh_0=\phi^{-1}X
\end{equation}
is the unit normal to the slices. We finally consider an $(n+1)$-dimensional Riemannian manifold $(\overline M,\overline g)$ endowed with a globally defined Killing field $\overline X$. We assume that there exists a closed subset $\overline K\subset \overline M$ such that $\overline M-\overline K$ is {\em isometric} to our warped product model $(\overline E\times I,\overline h)$, with   $\overline X$ corresponding to $X$ under the identification given by the isometry.

\begin{remark}\label{case}{\rm
An important special case of the above construction takes place when $(\overline M,\overline g)$ is {\em globally} a warped product. In this case, $(E,h)$ extends to a complete asymptotically flat manifold, still denoted $(E,h)$, and $\overline M=E\times I$ with $\overline g=h+{\phi}^2dt^2$, where now $\phi$ is a  positive extension of our function previously  defined only  on the end.
Here, $\overline X=\partial/\partial t$ everywhere and $\overline M$ is foliated by the totally geodesic slices $\overline E_{t}= E\times \{t\}$, $t\in I$, each of which has a globally defined unit normal vector field, namely,  $\eh_0=\phi^{-1}\overline X$. Notice that in this setting we can extend $R_h$ to $\overline M$ and ${\rm Ric}_{h}$ to act on arbitrary horizontal vectors in the obvious manner.}
\end{remark}

\begin{definition}\label{asympdef}
Let $(\overline M,\overline g)$ be as above. Then
a complete, isometrically immersed manifold $(M,g)\looparrowright (\overline M,\overline g)$, possibly with a compact inner boundary,   is {\em asymptotically flat} if $R_g$ is integrable and  there exists a compact subset $K\subset M$ such that $E_M=M-K$, the end of $M$, can be written as a vertical graph over some slice $E\hookrightarrow \overline M-\overline K$ associated to a smooth function $f:E\to\mathbb R$ such that the following asymptotic relations hold at infinity for some $\tau>(n-2)/2$:
\begin{equation}\label{asym}
\phi f_i(x)=O(|x|^{-\frac{\tau}{2}}), \quad \phi f_{ij}(x)+\phi_i f_{j}(x)=O(|x|^{-\frac{\tau}{2}-1}).
\end{equation}
Here, $f_i=\eh_i(f)=\partial f/\partial x_i$, etc.
\end{definition}

In view of (\ref{um}) and (\ref{mett1}) below, these decay conditions are tailored so that the mass of  $(M,g)$ is well defined and can be computed  using nonparametric coordinates at infinity according to the following slightly modified version of (\ref{mass}):
\begin{equation}\label{slight}
m_g=\lim_{r\to +\infty}c_n\int_{\Sigma_r}\left(\phi\left(g^j_{i,j}-g^j_{j,i}
\right)-\phi^jg_{ij}+\phi_ig^j_j\right)\nu^id\Sigma_r;
\end{equation}
in this regard, see \cite{M}.

We assume that $M$ is {\em two-sided} in the sense that it carries a globally defined  unit normal $N$, which we choose so that $N=\eh_0$ at infinity.
This allows us to consider the {\em angle function} $\Theta_{\overline X}=\langle \overline X,N\rangle:M\to\mathbb R$ associated to $\overline X$.
We also remark that the corresponding  orientation on $M$ induces, in the standard manner, an orientation on any inner boundary $\Gamma$. If we further assume  that $\Gamma$ lies on some (connected) hypersurface $P\hookrightarrow \overline M$ and there it bounds a distinguished compact domain $\Omega$, then a choice of a unit vector field, say $\mu$, pointing inward with respect to $\Omega$, defines a unique orientation on $P$.
With these preliminaries at hand, the following theorem computes, in the presence of an inner boundary $\Gamma$, which for simplicity we assume connected, the mass $m_g$ of $(M,g)$ in terms of the mass $m_h$ of $(E,h)$ and integrals over $\Gamma $ and $M$ involving geometric quantities.
In the following, we denote by ${\rm Ric}_{\overline g}$ (respectively, $R_{\overline g}$) the Ricci tensor (respectively, the scalar curvature) of $(\overline M,\overline g)$.

\begin{theorem}\label{mainwar}
Let $(M,g)\looparrowright (\overline M,\overline g)$ be a two-sided asymptotically flat hypersurface as above.
Let $\Gamma\hookrightarrow M$ be a smooth, compact inner boundary lying on some totally geodesic hypersurface $P\hookrightarrow \overline M$ and assume that, along $\Gamma$, $M$ is orthogonal to $P$. Then, if orientations are fixed as above,
\begin{eqnarray}\label{massformeuc}
m_g & = & m_h-c_n\int_{\Gamma}\langle \overline X ,\eta\rangle s_1(N) d\Gamma+\nonumber\\
 & & \qquad\quad +c_n\int_M\left(2S_2\Theta_{\overline X}+{\rm Ric}_{\overline g}(N,\overline X^T)\right) dM,
\end{eqnarray}
where $s_1(N)$ is the mean curvature of $\Gamma\hookrightarrow P$ with respect to $N$, $\eta$ is the exterior unit co-normal to $M$, $S_2$ is the $2$-mean curvature of $M$ (see (\ref{twomean}) below) and ${\overline X}^T$ is the tangential component of $\overline X$ along $M$.
\end{theorem}

\begin{remark}\label{ortho}{\rm We note that the assumption that $M$ meets $P$ orthogonally along $\Gamma$ implies that $\Gamma\hookrightarrow M$ is minimal, a geometric condition a horizon should necessarily satisfy. In fact, it corresponds to the requirement  that $|\nabla f|\to +\infty$ as $(x,f(x))$ approaches  $\Gamma$ in the graph case considered in \cite{L1}; see Remark \ref{gather} below for a clarification of this point.}
\end{remark}

\begin{remark}\label{several}{\rm
The proof of Theorem \ref{mainwar}  straightforwardly extends to the case in which  $M$ has finitely many ends. Thus, if we assume that each end, say $E^{(i)}_M$, of $M$ is asymptotically flat in the sense that it can be expressed as a graph over $E$ in terms of a function satisfying the decay conditions (\ref{asym}), then (\ref{massformeuc}) gets replaced by
\begin{eqnarray}\label{newww}
\sum_i\epsilon_im^{(i)}_g & = & m_h-c_n\int_{\Gamma}\langle\overline X,\eta\rangle s_1(N) d\Gamma+\nonumber\\
 & & \qquad +c_n\int_M\left(2S_2\Theta_{\overline X}+{\rm Ric}_{\overline g}(N,\overline X^T)\right)\, dM,
\end{eqnarray}
where the $m^{(i)}_g$'s are the masses attached to the ends in the usual manner and $\epsilon_i=\pm 1$ according to whether the unit normals to $E^{(i)}_M$ and $E$ agree or not at infinity.
Also, the obvious generalization of (\ref{newww})  holds in the presence of finitely many horizontal inner boundaries.}
\end{remark}

\begin{remark}\label{chn}{\rm
In recent years there has been much  interest in defining and computing mass-like invariants for  non-compact Riemannian manifolds whose geometry at infinity approaches some model geometry other than the Euclidean one. A notable example occurs in the so-called asymptotically hyperbolic case; see for instance \cite{CH}, \cite{CN}, \cite{H2} and \cite{M}. In this regard we note that the methods leading to (\ref{massformeuc}) are flexible enough do deal with this more general setting. In particular, integral formulae for the mass invariants of such hypersurfaces are also available, from which we are able to  draw interesting consequences like Positive Mass and Penrose-like inequalities. These results will be presented in a companion paper \cite{dLG}.
}
\end{remark}

Formula (\ref{massformeuc}) becomes specially interesting in the setting of Remark \ref{case}. In this case we will {\em always} assume that $P=E_{t_0}$, $t_0\in I$, i.e. $\Gamma$ lies on some horizontal slice. Also, given the distinguished compact domain $\Omega\subset E_{t_0}$ such that $\partial \Omega=\Gamma$, the orientation on $E_{t_0}$  corresponding to the choice of the inward pointing unit vector $\mu$ agrees with the one induced by $\eh_0$.  Under these conditions, it is easy to check that the possibilities $\eta=\pm \eh_0$ imply that $\langle \overline X,\eta\rangle=\pm  \phi$ and $N=\mp \mu$ along $\Gamma$, so that (\ref{massformeuc}) reduces to
\begin{eqnarray}\label{masswar}
m_g & = & m_h +c_n\int_{\Gamma} \phi s_1(\mu) d\Gamma+\nonumber\\
 & & \qquad\quad +c_n\int_M\left(2S_2\Theta_{\overline X}+{\rm Ric}_{\overline g}(N,\overline X^T)\right) dM,
\end{eqnarray}
where $s_1({\mu})$ is the mean curvature of $\Gamma \hookrightarrow E_{t_0}$ with respect to $\mu$.
This formula already comprises a number of interesting subcases that certainly deserve being discussed here, but for the sake of brevity we only mention one important consequence, whose complete justification is deferred to Proposition \ref{postp} below.

\begin{theorem}\label{lam2}
Assume in Theorem \ref{mainwar} that $(\overline M,\overline g)$ is a Riemannian product, i.e.  we are in the setting of Remark \ref{case} with $ \phi\equiv 1$. Then
\begin{eqnarray}\label{massformeuclam}
m_g & = & m_{h}+c_n\int_{\Gamma}s_1({\mu}) d\Gamma+\nonumber\\
 & & \qquad\quad +c_n\int_M\Theta_{\overline X}\left(R_g-R_{ h}+{\rm Ric}_{ h}(N^t,N^t)\right)\, dM,
\end{eqnarray}
where  $N^t$ is the horizontal component of $N$.
\end{theorem}

This has been proved in \cite{L2}
in the graph case (but with the boundary term missing). To explore (\ref{massformeuclam}) further, let us say that $(M,g)\hookrightarrow(\overline M,\overline g)$ as in Theorem \ref{lam2}
is a {\em quasi-graph} if $\Theta_{\overline X}\geq 0$ along $M$.
We also say that  $\Gamma\subset E_{t_0}$ is
{\em Alexsandrov-Fenchel} if $\Gamma$ is contained in a domain $D\subset E_{t_0}$ which is {\em isometric} to an Euclidean domain and  $\Gamma\subset D$ is convex (equivalently, the pair $\Gamma\subset D$ projects down to a pair $\Gamma_0\subset D_0$, where $D_0\subset E$, with the same properties, which means that this really amounts to an assumption on the geometry of  $(E,h)$). Then  the  Alexsandrov-Fenchel inequality applies to estimate from below the boundary integral in (\ref{massformeuclam}) in terms of the area $A$ of $\Gamma$ in the standard manner, which yields the following optimal {\em relative} Penrose-type inequality.

\begin{theorem}\label{lam3}
If, under the conditions of Theorem \ref{lam2}, $M$ is a quasi-graph and $\Gamma$ is Alexsandrov-Fenchel then
\begin{equation}\label{lam5}
m_g\geq m_h+\frac{1}{2}\left(\frac{A}{\omega_{n-1}}\right)^{\frac{n-2}{n-1}},
\end{equation}
whenever
\begin{equation}\label{hyp}
R_g\geq R_h-{\rm Ric}_h(N^t,N^t)
\end{equation}
outside of the zero set of $\Theta_{\overline X}$.
\end{theorem}

\begin{corollary}\label{classes}{
 Under the conditions of the theorem, assume that either $n\leq 7$ or  $E$ is spin and (\ref{hyp}) holds with $R_h\geq 0$. Then,
\begin{equation}\label{classes2}
m_g\geq \frac{1}{2}\left(\frac{A}{\omega_{n-1}}\right)^{\frac{n-2}{n-1}}.
\end{equation}
}
\end{corollary}

To check this, notice that $R_h\geq 0$ allows us to apply the Positive Mass Theorem \cite{SY} \cite{W} so that $m_h\geq 0$ and (\ref{classes2}) follows.

\begin{remark}\label{streng}{\rm
Theorem \ref{lam3} can be used to produce examples of asymptotically flat manifolds for which  the super-optimal inequality
\begin{equation}\label{superopt}
m_g\geq \left(\frac{A}{\omega_{n-1}}\right)^{\frac{n-2}{n-1}}
\end{equation}
holds.
For example, in dimension $n=3$, Bartnik \cite{Ba2} constructed examples of scalar-flat, asymptotically flat metrics containing a domain $D_0$ isometric to an Euclidean ball. Moreover, if combined with a result by Corvino \cite{C},  the metric can even be chosen to be isometric to the Schwarzschild solution in a neighborhood of infinity.  Clearly, the corollary applies to this class of metrics provided  $\Gamma\subset D$ and (\ref{hyp}) holds with $R_h=0$, thus yielding the super-optimal lower bound
\begin{equation}\label{super}
m_g\geq \sqrt{\frac{A}{4\pi}},
\end{equation}
since now both the horizon and the end of the quasi-graph each  contribute to its mass with the standard Penrose lower bound.
More generally, we can use a recent  gluing result due to Brendle, Marques and Neves \cite{BMN} to obtain similar examples in any dimension $n\geq 3$. The idea is to glue a \lq bowl\rq\,  metric in a spherical cap (so that the bottom of the bowl is flat) to the (exterior) Riemannian Schwarzschild solution along their common boundary, which is diffemorphic to $\mathbb S^{n-1}$. Since the boundary of the Schwarzschild solution is minimal, the conditions of Theorem 5 in \cite{BMN} are met, after possibly adjusting the size of the \lq bowl\rq, and in this way we obtain a manifold $(E,h)$ which contains  a flat region and is Schwarzschild at infinity, so we again can find quasi-graphs $(M,g)$ for which (\ref{superopt}) holds whenever
(\ref{hyp}) takes place. The negative part of the scalar curvature of the manifold $(E,h)$ so obtained  gets concentrated in a small neighborhood of the common boundary, but notice that we can arrange so that $R_h\geq 0$ everywhere if, before the gluing procedure, we slightly modify a small neighborhood of the boundary of the Schwarzschild solution in a rotationally invariant manner so as to render its scalar curvature positive (and small). That this kind of modification can be performed  is due to the fact that the mean curvature of the boundary of the \lq bowl\rq\, metric  is positive.}
\end{remark}

\begin{remark}\label{fs}{\rm
Another interesting example of a manifold to which we can attach  \lq bowl\rq\, metrics as in the previous remark appears in the conformally flat case. Let us consider a bounded domain $\Omega\subset \mathbb R^n$ and the manifold $(E',h')$, where $E'=\mathbb R^n-\Omega$ and $h'=u^{\frac{4}{n-2}}g_0$, with $g_0$ being the flat metric and $u>0$  a smooth function on $E'$ satisfying $\Delta_{g_0}u\leq 0$  everywhere and $u\to 1$ at infinity, i.e. $R_{h'}\geq 0$ and $h'$ is asymptotically flat; see \cite{BI}. We assume that $\partial E'$ is minimal in $(E',h')$ and convex as an Euclidean hypersurface. In this case  we can clearly attach a \lq bowl\rq\, type metric to $h'$ along $\partial E'$ so as to obtain a complete, asymptotically flat manifold $(E,h)$ whose scalar curvature is nonnegative except possibly in a small neighborhood of the common boundary. As in the previous remark, if before gluing we adjust the metric $h'$ so that its scalar curvature becomes  slightly positive in a neighborhood of the boundary, the condition $R_{h}\geq 0$ can be globally restored. Since the background manifold is spin, Corollary \ref{classes}  applies to a suitable asymptotically flat quasi-graph $(M,g)$ satisfying (\ref{hyp}).
But notice that, by using the recent lower bound for $m_{h'}$ obtained by Freire and Schwartz \cite{Sc} \cite{FS}, which of course turns into a lower bound for $m_h$, (\ref{classes2}) improves to
\begin{equation}\label{improv1}
m_g\geq 2\left(\frac{V(\Omega)}{\beta_n}\right)^{\frac{n-2}{n}}+
\frac{1}{2}\left(\frac{A}{\omega_{n-1}}\right)^{\frac{n-2}{n-1}},
\end{equation}
where $V(\Omega)$ (respectively, $\beta_n$) is the volume of $\Omega$ (respectively, of the unit ball). Alternatively, we may appeal to a recent result by Jauregui \cite{J} in order to get
\begin{equation}\label{improv2}
m_g\geq \left(\frac{|\partial E'|}{\omega_{n-1}}\right)^{\frac{n-2}{n-1}}+
\frac{1}{2}\left(\frac{A}{\omega_{n-1}}\right)^{\frac{n-2}{n-1}},
\end{equation}
where $|\partial E'|$ is the Euclidean area of $\partial E'$. Notice that, if combined with the classical isoperimetric inequality, (\ref{improv2}) yields a weakened version of (\ref{improv1}), with the factor $2$ in the first term on the right-hand side missing. We remark, however, that the lower bound  in \cite{J} applies to a much larger class of conformal metrics, and since our gluing procedure also works for such metrics, (\ref{improv2}) holds in this generality.}
\end{remark}

\begin{remark}\label{even}{\rm
Even in the rather special case $(E,h)=(\mathbb R^n, g_0)$ the corollary gives an optimal Penrose inequality in {\em any} dimension for quasi-graphs in $\mathbb R^{n+1}$ with $R_g\geq 0$ outside of the zero set of $\Theta_{\overline X}$ and whose inner boundary $\Gamma$ lies on a horizontal hyperplane $\Pi$, with $M$ meeting $\Pi$  orthogonally along $\Gamma$. This already strengthens a celebrated result verified for graphs in \cite{L1}. Notice that here the scalar curvature is even allowed to be {\em negative} somewhere in the zero set of $\Theta_{\overline X}$.}
\end{remark}

Another interesting consequence of Theorem \ref{mainwar} we mention here takes place when $(\overline M,\overline g)$ is Ricci-flat.

\begin{theorem}\label{ricci}
If $(\overline M,\overline g)$ is Ricci-flat, then
\begin{eqnarray}\label{massformeuc2}
m_g & = & m_h-c_n\int_{\Gamma}\langle \overline X,\eta\rangle s_1(N) d\Gamma+\nonumber\\
 & & \qquad\quad  +c_n\int_M\Theta_{\overline X}R_g dM.
\end{eqnarray}
\end{theorem}

In effect,  Ricci flatness  and (\ref{ger}) imply ${\rm Ric}_{\overline g}(N,X^T)=0$ and $2S_2=R_g$.
This yields a sort of positive mass inequality for certain hypersurfaces.

\begin{corollary}\label{classes3}{
Assume, under the conditions of the theorem, that $m_h\geq 0$, $\Gamma=\emptyset$ and $\Theta_{\overline X}\geq 0$ along $M$. Then, $m_g\geq 0$ whenever $R_g\geq 0$ outside of the zero set of $\Theta_X$.
}
\end{corollary}

We note that the conclusion holds if we merely assume that $\Theta_{\overline X}R_g\geq 0$ everywhere along $M$.
Also, it is not hard to check that (sub-optimal) versions of the Penrose inequality hold in the setting of the corollary under suitable assumptions
on a inner boundary $\Gamma$. For example, one might ask that $\Gamma$ lies on one of the slices of the standard foliation of  $\overline M-\overline K$. Since the slice is almost flat, Alexsandrov-Fenchel holds with an almost optimal constant if we assume that $\Gamma$ is \lq convex\rq\, in the obvious sense. Finally, we remark that an optimal Penrose inequality can be obtained if a neighborhood of $P$ is a Riemannian product, with $P$ as a slice, and $\Gamma$ is Alexsandrov-Fenchel.

The theorem also yields a scalar curvature rigidity result according to which one cannot, under proper conditions, deform the induced metric so that the scalar curvature increases at each point while keeping the geometry at infinity fixed. This is loosely related to results  surveyed in \cite{Bre}.

\begin{corollary}\label{rigid}{
Assume, under the conditions of the theorem, that $\Gamma=\emptyset$, $\Theta_{\overline X}>0$ along $M$, $f$ is constant in a neighborhood of infinity and $R_g\geq 0$. Then, $(M,g)$ is scalar-flat.}
\end{corollary}

Indeed, one has $m_g=m_h$, which gives
$$
\int_M\Theta_{\overline X}R_g dM=0.
$$

\begin{remark}\label{streng2}
We would like to point out that not every asymptotically flat manifold can be isometrically immersed in $(\mathbb{R}^{n+1},g_0)$ as an asymptotically flat hypersurface (in the sense of Definition \ref{asympdef}). This is certainly the case of the three-dimensional Bartnik-Corvino's examples in Remark \ref{streng}.
Indeed, if this were the case, the scalar-flatness of these manifolds would imply, by Lam's formula (\ref{lam}), that the mass of each of them is zero, a contradiction since they agree with a Schwarzschild solution at infinity.
Similar remarks also hold for immersions in the Ricci-flat manifolds of Theorem \ref{ricci}. This indicates the limitations of our methods and shows that our mass formulae can be seen as geometric obstructions to realizing certain asymptotically flat manifolds as asymptotically flat hypersurfaces in the ambient manifolds we consider.
\end{remark}

This paper is organized as follows. In Section \ref{geograph} we compute the shape operator of an asymptotically flat  graph in $(\overline M,\overline g)$. The result, presented in Proposition \ref{data} below, looks somewhat intractable at first sight, but we show in Proposition \ref{key} that a remarkable cancelation takes place if one evaluates the corresponding Newton tensor on the tangential component of the vertical Killing field. In Section \ref{proof} we then show, after integrating the flux formula by parts, that the above mentioned simplified expression, when restricted to large coordinate spheres, can be identified to the field of $1$-forms defining the ADM mass after discarding, by means of a careful analysis, higher order terms that vanish at infinity after integration over these spheres.
This proves Theorem \ref{mainwar} in case $\Gamma=\emptyset$ and an extra piece of argument is then presented to handle the general case. In Section \ref{quasilocal} we discuss a notion of quasi-local mass in the class of hypersurfaces we are considering, showing in particular that, under suitable conditions, this quasi-local mass is nonnegative, monotone and converges to the ADM mass. Finally, in Section \ref{further},  we briefly discuss further generalizations of our  results, including the case of (Riemannian or Lorentzian) manifolds carrying a conformal Killing field and whose geometry at infinity approaches  our warped product model in a suitable sense.

\vspace{0.1cm}
\noindent
{\bf Acknowledgements.} The authors would like to thank F. Marques for many helpful suggestions and, in particular, for pointing out the gluing result in \cite{BMN}. They are also indebted to F. Schwartz for valuable comments on an earlier version of this paper.

\section{The geometry of graphs in warped products}\label{geograph}

If $(M,g)\looparrowright(\overline M,\overline g)$ is a two-sided asymptotically flat  hypersurface as in Definition \ref{asympdef} and $\overline \nabla$ is the Riemannian connection of $(\overline M,\overline g)$, let us denote by  $B=-\overline \nabla N$ the shape operator of $M$ with respect to its unit normal vector $N$ and by  $k_1,\ldots,k_n$  the eigenvalues of $B$ with respect to $g$ (the principal curvatures). Define
\begin{equation}\label{onemean}
S_1=\sum_ik_i
\end{equation}
and
\begin{equation}\label{twomean}
S_2=\sum_{i<j}k_ik_j.
\end{equation}
These are respectively the {\em  mean curvature} and the $2$-{\em mean curvature} of $M$. Notice that from Gauss equation we have
\begin{equation}\label{ger}
R_g=R_{\overline g}-2{\rm Ric}_{\overline g}(N,N)+2S_2.
\end{equation}
Also, we define  the {\em Newton tensor} by
\begin{equation}\label{zero}
G=S_1I-B,
\end{equation}
where $I$ is the identity map.

Later on we will need the expressions of some of these invariants along the end  $E_M$ of $M$ which, by Definition \ref{asympdef}, is a graph over the end  $E$. To achieve  this we start by noticing that if, as before,  $x=(x_1,\cdots,x_n)$ are asymptotically flat coordinates in $E$, then the tangent frame
\begin{equation}\label{tangent}
\eh_i=\frac{\partial}{\partial x_i}, \quad i=1,\cdots,n,
\end{equation}
can be extended to a frame $\{\eh_{\alpha}\}_{\alpha=0}^n$ in $\overline E=E\times\mathbb R$ with $\langle\eh_0,\eh_i\rangle=0$, where $\eh_0$ is defined by ({\ref{framevert}}). The following proposition describes the structure equations associated to such a frame.

\begin{proposition}\label{struct}
One has
\begin{equation}\label{structw}
\overline\nabla_{\eh_i}\eh_0=0,\quad \overline\nabla_{\eh_0}\eh_i=\phi^{-1}\phi_i\eh_0,\quad
\overline\nabla_{\eh_0}\eh_0=-\phi^{-1}\nabla \phi,
\end{equation}
where $\nabla$ is the gradient operator of $(E,h)$ and $\phi_i=\eh_i(\phi)$.
\end{proposition}

\begin{proof}
The first equation in (\ref{structw}) follows from the fact that the slices are totally geodesic, as already remarked.
From this we get
\begin{eqnarray*}
\overline\nabla_{\e_0}\e_i=\overline\nabla_{\e_i}\e_0+[\e_0,\e_i] & = & \left[\phi^{-1}\frac{\partial}{\partial t},\e_i\right]\\
& = &
\phi^{-1}\left[\frac{\partial}{\partial t},\e_i\right]-\e_i(\phi^{-1})\frac{\partial}{\partial t}\\
& = & -\e_i(\phi^{-1})\frac{\partial}{\partial t},
\end{eqnarray*}
and the second equation follows. Finally,
$$
\overline\nabla_{\e_0}\e_0=\phi^{-2}\overline\nabla_{\partial/\partial t}
\frac{\partial}{\partial t},
$$
and, since $\partial/\partial t$ is Killing, this implies
$$
\langle \overline\nabla_{\e_0}\e_0,\e_0\rangle=\phi^{-3}
\left\langle\overline\nabla_{\partial/\partial t}\frac{\partial}{\partial t},
\frac{\partial t}{\partial t}\right\rangle=0,
$$
where $\langle\,,\,\rangle$ is the inner product associated to $\overline g$. Thus,
$$
\overline\nabla_{\e_0}\e_0=\gamma^l\e_l,
$$
with
\begin{eqnarray*}
\gamma^l & = & h^{lm}\langle\overline\nabla_{\e_0}\e_0,\e_m\rangle \\
& = & -\phi^{-2}
h^{lm}\langle\overline\nabla_{\e_m}\frac{\partial}{\partial t},
\frac{\partial}{\partial t}\rangle\\
& = & -\frac{\phi^{-2}}{2}{h^{lm}}
\e_m\left(\phi^2\right)\\
& = & -\phi^{-1}h^{lm}\phi_m,
\end{eqnarray*}
as desired.
\end{proof}

Let us now write
$$
E_M=\left\{(x,f(x));x\in E\right\}\subset\overline M,
$$
as the graph associated to a smooth function $f:E\to\mathbb R$ as in Definition \ref{asympdef}.
In terms of the frame in Proposition \ref{struct}, $TE_M$ is spanned by
\begin{equation}\label{frame}
Z_i=f_i\frac{\partial}{\partial t}+\e_i=\phi f_i\e_0+\e_i,\quad i=1,\cdots,n,
\end{equation}
where $f_i=\eh_i(f)$, and we  choose
\begin{equation}\label{normal}
N=\frac{1}{W}\left(\e_0-\phi \nabla f\right),
\end{equation}
where
\begin{equation}\label{asymw}
W=\sqrt{1+\phi^2|\nabla f|^2_h}=1+O(|x|^{-\tau}),
\end{equation}
as the unit normal to $E_M$. Notice that this is consistent with our global choice of unit normal to $M$, which is dictated by the condition $N=\eh_0$ at infinity.
Also, the induced metric on $E_M$ is
\begin{equation}\label{mett1}
g_{ij}=\langle Z_i,Z_j\rangle=h_{ij}+\phi^2f_if_j,
\end{equation}
where $h_{ij}=\langle\eh_i,\eh_j\rangle$ is the metric on $E$, and
the inverse metric is
\begin{equation}\label{mett2}
g^{ij}=h^{ij}-\frac{\phi^2}{W^2}f^if^j,
\end{equation}
where here and everywhere else in the paper indexes are raised and lowered using  $h$.

\begin{proposition}\label{data}
The coefficients of the shape operator $B$ of the graph $E_M$  with respect to the frame (\ref{frame}) are
\begin{eqnarray}\label{shape}
WB_{j}^i & = & h^{ik}\left(\phi f_{kj}+\phi_kf_j+\phi_jf_k+\phi^2f_kf_j\phi^mf_m\right)
                 -\nonumber\\
&  & \quad - \frac{\phi^2}{W^2}f^if^k\left(\phi f_{kj}+\phi_kf_j+\phi_jf_k+\phi^2f_kf_j\phi^mf_m\right)
\label{data2}.
\end{eqnarray}
\end{proposition}

\begin{proof}
We shall use (\ref{structw}) and  start by computing the coefficients
$$
\alpha_{jk}=\langle \overline\nabla_{Z_j}Z_k,N\rangle
$$
of the second fundamental form $\alpha$ of $E_M$. Since $\overline\nabla_{\e_j}\e_0=0$, a direct computation gives
$$
\overline\nabla_{Z_j}Z_k=\phi f_j\e_0(\phi f_k)\e_0+\phi^2f_jf_k\overline\nabla_{\e_0}\e_0+
\phi f_j\overline\nabla_{\e_0}\e_k+
\e_j(\phi f_k)\e_0+\overline\nabla_{\e_j}\e_k.
$$
But notice that $\e_0(\phi f_k)=\phi^{-1}\partial_t(\phi f_k)=0$. Moreover,
$$
\overline\nabla_{\e_j}\e_k=\nabla_{\e_j}\e_k+\beta^{(t)}(\e_j,\e_k),
$$
where $\nabla$ and $\beta^{(t)}$ are the connection and second fundamental form of a slice $E_t\subset \overline M$. Since $\beta^{(t)}=0$ and we may assume that $\nabla_{\e_j}\e_k=0$ at the point where we are doing the computation, it follows that $\overline\nabla_{\e_j}\e_k=0$. Thus,
$$
\overline\nabla_{Z_j}Z_k=\phi^2f_jf_k\overline\nabla_{\e_0}\e_0+
\phi f_j\overline\nabla_{\e_0}\e_k+
\e_j(\phi f_k)\e_0,
$$
and from (\ref{structw}) and (\ref{normal}) we easily get
\begin{equation}\label{braykhuri}
\alpha_{jk}=\frac{1}{W}\left(\phi f_{jk}+\phi_jf_k+\phi_kf_j+\phi^2f_jf_k\phi^mf_m\right).
\end{equation}
In view of (\ref{mett2}), the expression (\ref{shape}) for the shape operator $B_{k}^i=g^{ij}\alpha_{jk}$  follows readily.
\end{proof}

\begin{remark}\label{braykhurirem}{\rm
A computation leading to the Lorentzian analogue of (\ref{braykhuri}) is presented in \cite{BK}.
}
\end{remark}

The following proposition is a key ingredient in our proof of Theorem \ref{mainwar} as it shows that the specific combination of extrinsic data yielding the Newton tensor  of a graph simplifies considerably after evaluation on the tangential component of the vertical Killing field.

\begin{proposition}\label{key}
Let $E_M\subset\overline M$ be an asymptotically flat graph as above, $G$ its Newton tensor and $X^T$ the tangential component of $X=\partial/\partial t$ along $E_M$. Then, with respect to the frame (\ref{frame}), the coefficients of $GX^T$ are given by
\begin{equation}\label{kx1}
(GX^T)^i=(GX^T)^i_{(1)}+(GX^T)^i_{(2)},
\end{equation}
where
\begin{equation}\label{kx2}
(GX^T)^i_{(1)}= \frac{\phi^2}{W^3}\left(\phi f_{kj}+\phi_kf_j+\phi_jf_k\right)\left(h^{jk}f^i-h^{ik}f^j\right)=O(|x|^{-\tau-1}),
\end{equation}
and
\begin{equation}\label{kx4}
(GX^T)^i_{(2)}=\frac{\phi^4}{W^3}\phi^mf_mf_kf_j\left(h^{jk}f^i-
 h^{ik}f^j\right)=O(|x|^{-2\tau-1}).
\end{equation}
\end{proposition}

\begin{proof}
From (\ref{zero}) we have
\begin{equation}\label{form}
(GX^T)^i=B_{j}^j(X^T)^i-B^i_j(X^T)^j,
\end{equation}
where
\begin{equation}\label{comp2}
X^T=(X^T)^iZ_i=(X^T)^i\eh_i+\phi f_i(X^T)^i\eh_0
\end{equation}
by (\ref{frame}).
Let us rewrite (\ref{shape}) as
$$
B_{j}^i=\sum_{s'=1}^8{B_{j}^i}_{(s')},
$$
where $W{B_{j}^i}_{(1)}=\phi h^{ik}f_{kj}$, $W{B_{j}^i}_{(2)}=h^{ik}\phi_kf_j$, etc.
Now, since $\langle X,N\rangle=\phi/W$,
$$
X^T=X-\frac{\phi}{W}N=\frac{\phi^3 |\nabla f|^2_h}{W^2}\e_0+\frac{\phi^2}{W^2}f^i\eh_i,
$$
and comparing with (\ref{comp2}),
$$
(X^T)^i=\frac{\phi^2}{W^2}f^i.
$$
It is now straightforward to check that
$$
{B_{j}^j}_{(s)}(X^T)^i={B_{j}^i}_{(s)}(X^T)^j,\quad s'\ge 5,
$$
that is, half the terms in (\ref{form}) cancel out
and (\ref{kx1}) follows easily. The decay rates follow from the corresponding ones in Definition \ref{asympdef} and the fact that both $\phi$ and $h$ are uniformly bounded at infinity.
\end{proof}

\section{The proof of Theorem \ref{mainwar}}\label{proof}

In this section we present the proof of Theorem \ref{mainwar}. As remarked in the Introduction, the starting point is the flux-type formula
\begin{equation}\label{basic}
{\rm div}_gG\overline X^T=2S_2\Theta_{\overline X} + {\rm Ric}_{\overline g}(N,\overline X^T),
\end{equation}
where $(M,g)\looparrowright (\overline M,\overline g)$ is a two-sided asymptotically flat hypersurface as in Definition \ref{asympdef}, $G$ is its Newton tensor and $\overline X^T$ is the tangential component of the Killing field $\overline X$ that agrees  with $X=\partial/\partial t$ on $\overline M-\overline K$.
In this generality, (\ref{basic}) has been first obtained in \cite{ABC} in the Lorentzian setting. The Riemannian version can be found in \cite{AdLM}; see their Lemma 3.1 and equation (8.4) with $r=1$, but be aware of their choice for the sign of the curvature tensor.

Before proceeding with the proof, let us consider the  special form of  (\ref{basic}) leading to Theorem \ref{lam2} in the Introduction.

\begin{proposition}\label{postp}
Under the conditions of Theorem \ref{lam2},
(\ref{basic}) reduces to
\begin{equation}\label{basic2}
{\rm div}_gG\overline X^T=\Theta_{\overline X}\left(R_g- R_h +{\rm Ric}_h(N^t,N^t)\right).
\end{equation}
\end{proposition}

\begin{proof}
The product structure gives ${R}_{\overline g}=R_h$ and ${\rm Ric}_{\overline g}(\eh_0,Y)=0$ for any $Y$, so if $N=N^t+N^n$, with $N^n$ proportional to $\eh_0$, we have
\begin{equation}\label{then}
{\rm Ric}_{\overline g}(N,N)  =  {\rm Ric}_{\overline g}(N^t,N^t) =  {\rm Ric}_h(N^t,N^t),
\end{equation}
where, in the last step, we have used that the slices are totally geodesic.
On the other hand, since $\overline X^T=\overline X-\Theta_{\overline X}N$,
\begin{eqnarray*}
{\rm Ric}_{\overline g}(N,\overline X^T) & = & \phi{\rm Ric}_{\overline g}(N,\eh_0)-
           \Theta_{\overline X}{\rm Ric}_{\overline g}(N,N)\\
           & = & -\Theta_{\overline X}{\rm Ric}_h(N^t,N^t),
\end{eqnarray*}
and (\ref{basic2})  follows from (\ref{ger}), (\ref{basic}) and (\ref{then}).
\end{proof}

Returning to the proof of Theorem \ref{mainwar}, we first consider the case $\Gamma=\emptyset$, i.e. no inner boundary is present. We take a large coordinate sphere $\Sigma_r\subset E$  and set $\sigma_r=f(\Sigma_r)$,
where $f$
describes $E_M$ as a graph over $E$.
We denote by $M_r$ the compact region of $M$ inside $\sigma_r$ and by $\vartheta$ (respect. $\nu$) the outward unit normal to $\sigma_r$ (respect. $\Sigma_r$). Thus, integrating (\ref{basic}) over $M$ and using the divergence theorem, we get
\begin{eqnarray}\label{div}
 \int_M\left(2S_2\Theta_{\overline X}+{\rm Ric}_{\overline g}(N,\overline X^T)\right)dM & = & \lim_{r\to \infty}\int_{\sigma_r}\langle G\overline X^T,\vartheta\rangle \,d\sigma_r\nonumber\\
& = &  \lim_{r\to \infty}\int_{\Sigma_r}g_{im}(G\overline X^T)^i\nu^m\,d\Sigma_r,
\end{eqnarray}
where we used that at infinity we may replace $\vartheta^md\sigma_r$ by $\nu^md\Sigma_r$.
Thus, comparing with (\ref{massformeuc}) we  are left with the task of relating the right-hand side above, which manifestly depends on the extrinsic geometry of the end, to the intrinsically defined relative mass $m_g-m_h$.

By (\ref{mett1}) and Proposition \ref{key},
\begin{equation}\label{disc}
g_{im}(G\overline X^T)^i\nu^m = h_{im}\nu^m\sum_{s=1}^2(G\overline X^T)^i_{(s)} +\phi^2f_if_m\nu^m\sum_{s=1}^2(G\overline X^T)^i_{(s)},
\end{equation}
but notice that (\ref{kx2}) implies
$$
\phi^2f_if_m\nu^m(G\overline X^T)^i_{(1)}=O(|x|^{-2\tau-1}),
$$
and since the area of coordinate spheres in $E$ grows as $|x|^{n-1}$, this term vanishes at infinity after integration because it becomes $O(|x|^{-2\tau+n-2})$ there.
Similarly, by (\ref{kx4}),
$$
\phi^2f_if_m(G\overline X^T)^i_{(2)}\nu^m=O(|x|^{-3\tau-1}),
$$
and we obtain
\begin{equation}\label{disc2}
g_{im}(G\overline X^T)^i\nu^m \approx h_{im}\nu^m\sum_{s=1}^2(G\overline X^T)^i_{(s)},
\end{equation}
where $\approx$ means precisely that we are discarding terms that vanish at infinity after integration.

In order to get rid of further terms in (\ref{disc2}) we first note that (\ref{decay1}) implies
$$
h_{im}(G\overline X^T)^i_{(s)}\nu^m\approx (G\overline X^T)^i_{(s)}\nu_i,\quad 1\leq s\leq 2.
$$
But by (\ref{kx4}),
$$
(G\overline X^T)^i_{(2)}\nu_i\approx 0,
$$
and using that $h^{ij}=\delta^{ij}+O(|x|^{-\tau})$, which follows from (\ref{decay1}), we  obtain
$$
(G\overline X^T)^i_{(1)}\nu_i\approx I(\phi)^i\nu_i,
$$
where
$$
I(\phi)^i  =
\frac{\phi^2}{W^3}\left(\phi\left(f_{j}^jf^i-f_{j}^if^j\right)+
  \phi^jf_jf^i-\phi^if_jf^j\right).
$$
Thus, we find that (\ref{div}) can be rewritten as
\begin{equation}\label{div2}
\int_M\left(2S_2\Theta_{\overline X}+{\rm Ric}_{\overline g}(N,\overline X^T)\right)dM=\lim_{r\to +\infty}
\int_{\Sigma_r}I(\phi)^i\nu_id\Sigma_r.
\end{equation}

To relate the integrand in the right-hand side of (\ref{div2})  to  the $1$-form defining the mass we now observe that
(\ref{slight}), (\ref{asymw}) and (\ref{mett1}) lead to
\begin{equation}\label{quasi}
m_g=m_h+ \lim_{r\to +\infty}c_n\int_{\Sigma_r}J(\phi)_i\nu^id\Sigma_r,
\end{equation}
where
$$
J(\phi)_i=\frac{1}{W^3}\left(\phi\left(e^j_{i,j}-e^j_{j,i}\right)-\phi^je_{ij}+
\phi_ie^j_j\right)
$$
and $e_{ij}=g_{ij}-h_{ij}=\phi^2f_if_j$. Now a straightforward computation gives $J(\phi)_i\nu^i=I(\phi)^i\nu_i$,
so that (\ref{quasi}) becomes
$$
\lim_{r\to \infty}c_n\int_{\Sigma_r}I(\phi)^i\nu_i\,d\Sigma_r =m_g-m_h,
$$
which together with (\ref{div2}) completes the proof of Theorem \ref{mainwar} in case $\Gamma=\emptyset$.

In the presence of $\Gamma$, the  extra boundary integral
$$
-\int_{\Gamma}\langle G\overline X^T,\eta\rangle d\Gamma,
$$
where $\eta$ is the outward unit co-normal to $\Gamma$, pops out in the left-hand side of (\ref{div2}). To properly handle this we use our orthogonality assumption to expand, in terms of a local orthonormal basis
$\{\tilde\eh_l \}_{l=1}^{n-1}$ of $T\Gamma$,
$$
\overline X^T  = \langle \overline X,\eta\rangle\eta+\sum_l\langle \overline X,\tilde \eh_l\rangle\tilde \eh_l,
$$
so that
\begin{eqnarray*}
\langle G\overline X^T,\eta\rangle &  = &  \langle \overline X,\eta\rangle\langle G\eta,\eta\rangle+\sum_l\langle \overline X,\tilde \eh_l\rangle\langle G\tilde \eh_l,\eta\rangle\\
  & = & \langle \overline X,\eta\rangle(S_1-\langle B\eta,\eta\rangle)+\sum_l\langle \overline X,\tilde \eh_l\rangle\langle G\tilde \eh_l,\eta\rangle.
\end{eqnarray*}
But
$$
\langle G\tilde\eh_l,\eta \rangle=-\langle B\tilde\eh_l,\eta\rangle=\langle \overline\nabla_{\tilde \eh_l}N,\eta\rangle= -\langle N,\overline\nabla_{\tilde \eh_l}\eta\rangle,
$$
and this vanishes due to the assumptions that $\eta=\pm\xi$ along $\Gamma$ and that $P$ is totally geodesic. In particular, $\eta$ is a principal direction of $B$ with $\langle B\eta,\eta\rangle$  being the corresponding principal curvature and hence $S_1-\langle B\eta,\eta\rangle=s_1(N)$, which completes the proof of  Theorem \ref{mainwar}.

\begin{remark}\label{novel}
{\rm A rather informal, but highly suggestive, way of concisely expressing the mass is
$$
m_g=c_n\int_{\Sigma_\infty} \langle\Upsilon_{ADM},\nu_{\infty}\rangle d\Sigma_{\infty},
$$
where $\Upsilon_{ADM}={{\rm div}_{g_0}g-d{\rm tr}_{g_0}g}$, $\Sigma_\infty=\lim_{r\to +\infty} \Sigma_r$ is the sphere at infinity and $\nu_\infty$ is its unit normal. In words, $m_g$ is simply the (properly normalized) total flux of the (co-)vector $\Upsilon_{ADM}$ over $\Sigma_{\infty}$. In this regard, the computation leading to the proof of Theorem \ref{mainwar} essentially amounts to checking that
$$
\int_{\Sigma_\infty} \langle\Upsilon_{ADM},\nu_{\infty}\rangle d\Sigma_{\infty}=\int_{\Sigma_\infty} \langle GX^T,\nu_{\infty}\rangle d\Sigma_{\infty},
$$
i.e. $\Upsilon_{ADM}$ and $GX^T$ have the {\em same} total flux over $\Sigma_\infty$. In fact, we have proved that $\Upsilon_{ADM}=GX^T+Y$, where $Y$, which corresponds to terms vanishing at infinity after integration, has a null total flux.}
\end{remark}

\begin{remark}\label{gather}{\rm
In order to justify the claim in Remark \ref{ortho} let us assume that $M$ meets $P$ in a not necessarily orthogonal manner. Thus, retaining the notation above, one has for each $l=1,\cdots,n-1$,
$$
\overline \nabla_{\tilde\eh_{l}}\tilde\eh_{l}=\sum_m\langle \overline\nabla_{\tilde\eh_{l}}\tilde\eh_{l},\tilde\eh_m\rangle\tilde\eh_m+
\langle\overline\nabla_{\tilde\eh_{l}}\tilde\eh_{l},N\rangle N+
\langle\overline\nabla_{\tilde\eh_{l}}\tilde\eh_{l},\eta\rangle \eta
$$
and
$$
\overline \nabla_{\tilde\eh_{l}}\tilde\eh_{l}=\sum_m\langle \overline\nabla_{\tilde\eh_{l}}\tilde\eh_{l},\tilde\eh_m\rangle\tilde\eh_m+
\langle\overline\nabla_{\tilde\eh_{l}}\tilde\eh_{l},\mu\rangle \mu+
\langle\overline\nabla_{\tilde\eh_{l}}\tilde\eh_{l},\xi\rangle \xi,
$$
where $\xi$ is the unit normal to $P$.
Taking the inner product of both equations with $\eta$, summing over $l$ and using that $P$ is totally geodesic we obtain
\begin{equation}\label{orthono}
s_1(\eta)=s_1(\mu)\langle \mu, \eta\rangle,
\end{equation}
where $s_1(\eta)$ is the mean curvature of $\Gamma\hookrightarrow M$. The claim follows.
}
\end{remark}

\section{A generalization of Lam's quasi-local mass}
\label{quasilocal}

As evidenced by the various Positive Mass and Penrose inequalities available in the literature, the ADM mass of an asymptotically flat manifold provides a rather satisfactory description of the energy content of an isolated gravitational system in General Relativity. It is highly desirable, however, to develop a notion of mass at the quasi-local level, i.e. for finitely extended regions in space. Such a notion should be expressed solely in terms of the boundary data and is required to meet a list of natural properties such as positivity, strict monotonicity, etc. In this regard we should mention that several proposals have been considered so far but a completely satisfactory solution to the problem remains elusive; see \cite{Ba3} and \cite{Sz} for surveys on this subject. The purpose of this section is to present a new notion of quasi-local mass for certain bounded domains in asymptotically flat hypersurfaces. As checked below, this quasi-local mass is nonnegative, monotone and ADM convergent under suitable conditions but, as already pointed out in the Introduction, our method fails to detect whether it is positive or strictly monotone in general, since this is essentially a rigidity issue.

In \cite{L2} it is defined a notion of {\em quasi-local mass} for bounded domains in an asymptotically flat  graph $M\hookrightarrow\mathbb R^{n+1}$ whose boundary $\Gamma$ is not necessarily a horizon, which means that $\Gamma$ is contained in a horizontal hyperplane $\Pi$ but $M$ is not orthogonal to $\Pi$ along $\Gamma$.
We will now show how this concept can be extended, in the presence of the Killing field $\overline X$, to certain bounded domains $\Dc\subset M$, where  $M\looparrowright\overline M$ is an asymptotically flat hypersurface . We will make two basic assumptions here. First, we assume that $\Gamma=\partial \Dc$ is the intersection of $M$ with a totally geodesic, embedded hypersurface $P\hookrightarrow\overline M$, as in our general setup. Second, we require that the Killing field $\overline X$ is normal to $P$ along $\Gamma$. We will then say that $(P,\Gamma,\Dc)$ is an {\em admissible configuration}. Notice that we do {\em not} assume that  $M$ is orthogonal to $P$ along $\Gamma$, so that $\Gamma\hookrightarrow M$ is not necessarily minimal; see Remarks \ref{ortho} and \ref{gather}. In any case, under these conditions  we have
$$
\overline X=\langle \overline X,\xi\rangle\xi=\langle \overline X,\eta\rangle\eta+
\langle\overline X,N\rangle N,
$$
where, as usual, $\eta$ is the outward unit co-normal to $M-\Dc$,
so that
\begin{equation}\label{notnece}
\langle G\eta,\overline X\rangle=\langle\overline X,\xi\rangle\langle\xi,\eta\rangle\langle G\eta,\eta\rangle.
\end{equation}
But formula 6.4 in \cite{AdLM} with $r=1$ says that $\langle G\eta,\eta\rangle=-s_1({\mu})\langle\xi,\eta\rangle$ and using that
$\langle \xi,\eta\rangle^2+\langle \eta,\mu\rangle^2=1$ and (\ref{orthono}) with $s_1(\mu)\neq 0$, we finally obtain
\begin{equation}\label{notnece2}
\langle G\eta,\overline X\rangle=\langle \overline X,\xi\rangle\left(
\frac{1}{s_1(\mu)}\left(s_1(\eta)^2-s_1(\mu)^2\right)\right).
\end{equation}
This motivates the following extension of the notion of quasi-local mass introduced in \cite{L2}.

\begin{definition}\label{qlmass}
Under the conditions above, the {\em quasi-local mass} of $\Dc$ is
\begin{equation}\label{qlmassf}
m_{QL}(\Dc)=c_n\int_{\Gamma}\langle \overline X,\xi\rangle\left(
\frac{1}{s_1(\mu)}\left(s_1(\mu)^2-s_1(\eta)^2\right)\right)d\Gamma.
\end{equation}
\end{definition}

Notice that, by (\ref{notnece2}),
$$
m_{QL}(\Dc)=-\int_{\Gamma}\langle G\overline X^T,\eta\rangle d\Gamma,
$$
so that, by (\ref{basic}) and the divergence theorem,
\begin{equation}\label{mql1}
m_{QL}(\Dc)=c_n\int_\Dc\left(2S_2\Theta_{\overline X}+{\rm Ric}_{\overline g}(N,\overline X^t)\right)dM.
\end{equation}
Also, the computation leading to the proof of Theorem \ref{mainwar} yields
\begin{equation}\label{mql2}
m_g=m_h+m_{QL}(\Dc)+c_n\int_{M-\Dc}\left(2S_2\Theta_{\overline X}+{\rm Ric}_{\overline g}(N,\overline X^t)\right)dM.
\end{equation}
More precisely, (\ref{mql2}) follows from (\ref{mql1}) and (\ref{massformeuc}) in case no horizon is present.

Here is the first consequence of our computation.

\begin{theorem}\label{conv} (Convergence to the ADM mass) If $(P_k,\Gamma_k,\Dc_k)$
is a sequence of admissible configurations with $\Dc_k$ exhausting $M$ as $k\to+\infty$ then
$$
\lim_{k\to+\infty}m_{QL}(\Dc_k)=m_g-m_h.
$$
\end{theorem}

\begin{proof}
This follows immediately from (\ref{mql2}) since
$$
\lim_{k\to +\infty} \int_{M-\Dc_k}\left(2S_2\Theta_{\overline X}+{\rm Ric}_{\overline g}(N,\overline X^t)\right)dM=0.
$$
\end{proof}

Further properties of $m_{QL}$ can be derived from (\ref{mql1}) under the conditions of Theorems \ref{lam2} and \ref{ricci}, as the following result shows.

\begin{theorem}\label{positi}
Assume that either: i) $(\overline M,\overline g)$ is a Riemannian product, $M$ is a quasi-graph and $R_g\geq R_h-{\rm Ric}_h(N^t,N^t)$ outside of the zero set of $\Theta_{\overline X}$ or ii) $(\overline M,\overline g)$ is Ricci-flat, $\Theta_{\overline X}\geq 0$ along $M$  and $R_g\geq 0$ outside of the zero set of $\Theta_{\overline X}$.
Then the following properties hold:
\begin{enumerate}
 \item (nonnegativity) If $(P,\Gamma,\Dc)$ is  admissible then $m_{QL}(D)\geq 0$;
  \item (monotonicity) If $(P_k,\Gamma_k,\Dc_k)$, $k=1,2$, are admissible configurations with $\Dc_1\subset \Dc_2$ then $m_{QL}(\Dc_1)\leq m_{QL}(\Dc_2)$.
\end{enumerate}
\end{theorem}

\begin{proof}
This follows from (\ref{mql1}) since in both cases the integrand in the right-hand side is non-negative.
\end{proof}

We remark that, in the case of Euclidean graphs, Theorems \ref{conv} and \ref{positi} have been proved in \cite{L2}.

\section{Further generalizations}\label{further}

In this section we briefly discuss a few generalizations of the results presented above.

\subsection{The case of conformal Killing fields}

The bulk of the argument leading to Theorem \ref{mainwar} actually involves a computation at infinity which  only uses suitable decay assumptions on the asymptotic geometry of the various manifolds involved. This suggests that a generalization of Theorem \ref{mainwar} should hold in cases where,  while maintaining the warped product structure at infinity, more flexibility is allowed on the geometry at finite scales. We briefly discuss here, in a rather sloppy style, one such possibility. Thus assume that $(\overline M_o,\overline g_o)$ is a Riemannian manifold for which there exists a closed subset $\overline K_o\subset \overline M_o$ with $\overline M_o-\overline K_o$ diffeomorphic to $\{(y,t)\in \mathbb R^n\times I;|y|>1\}$. Assume that $\overline M_o$ carries a conformal Killing field $\overline X_o$ and that, as one approaches infinity along $\overline M_o-\overline K_o$, the geometry of the triple $(\overline M_o,\overline g_o,\overline X_o)$ converges in a suitable sense to the geometry of our model $(\overline E,\overline h,X)$ for a function $\phi$ satisfying (\ref{um}). For example, a simple way of meeting these conditions is to require that at a neighborhood of spatial infinity $\overline X_0$ is Killing and the geometries are isometric indeed. Now, if we fix a copy of $E$ inside $(\overline M_o,\overline g_o)$, so that it has a well-defined mass $m_{h}$, let $(M_o,g_o)$ be an asymptotically flat hypersurface in $(\overline M_o,\overline g_o)$ which at infinity is a graph over (the fixed copy of) $E$ for some smooth function $f$  decaying as in (\ref{asym}). Thus, using the appropriate version of the flux formula (see \cite{AdLM}) and arguing  as above we will eventually find in  the case $\Gamma=\emptyset$ the following formula for the mass $m_{g_o}$ of $(M_o,g_o)$:
\begin{equation}\label{flex}
m_{g_o}  =  m_{h} + c_n\int_{M_o}\left(\lambda S_1+\left(2S_2\Theta_{\overline X_o}+{\rm Ric}_{{\overline g}_o}(N_o,\overline X_o^T)\right)\right) dM_o,
\end{equation}
where $\lambda$ is the conformality factor of $\overline X_o$,  $N_o$ is a suitably chosen unit normal to $M_o$ and the extrinsic invariants $S_1$ and $S_2$ now refer to $M_o$.
Thus, if $\overline M_o$ is Ricci-flat  and $M_o$ is minimal, we obtain the following analogue of (\ref{massformeuc2}):
$$
m_{g_o} =  m_{h}  +c_n\int_M\Theta_{\overline X_o}R_{g_o} dM_o.
$$
In particular, if $\Theta_{\overline X_o}\geq 0$ everywhere, we get $m_{g_o}\geq m_{h}$ provided $R_{g_o}\geq 0$ outside the zero set of $\Theta_{X_o}$.
Needless to say, under suitable assumptions, versions of (\ref{flex}) in the presence of an inner boundary can  be easily derived as well.

\subsection{The Lorentzian case}

Our main result, Theorem \ref{mainwar}, also admits a version in the Lorentzian case, which is obtained from our general Riemannian setup by \lq Wick rotation\rq. More precisely, $(\overline M,\overline g)$ now is a Lorentzian manifold which at \lq spatial\rq\, infinity agrees with the warped product model $(\overline E,\overline h)$, where $\overline E=E\times I$ and $\overline h=h-\phi^2 dt^2$. Additionally, we assume the existence of a time-like (conformal) Killing field
$\overline X$ agreeing with $\partial/\partial t$ in a neighborhood of infinity. The definition of asymptotically flat hypersurfaces $M\looparrowright \overline M$ is then the same as before, except that now $M$ is required to be space-like. Taking into account the results in \cite{ABC}, it is not hard to check that all the results above have counterparts in this Lorentzian setting.

\bibliographystyle{amsplain}

\begin{thebibliography}{999999}

\markboth{}{}



\bibitem{ABC} L. J. Al\'{\i}as, A. Brasil, Jr. and A. G. Colares, Integral formulae for spacelike hypersurfaces in conformally stationary spacetimes and applications. {\em Proc. Edinb. Math. Soc.} (2) 46 (2003), 465-488.


\bibitem{AdLM} L. J. Al\'{\i}as, J. H. S. de Lira and J. M. Malacarne,  Constant higher-order mean curvature hypersurfaces in Riemannian spaces, {\em J. Inst. Math. Jussieu} (2006), no. 4, 527-562.

\bibitem{ADM} R. Arnowitt, S. Deser and C. W. Misner,
Energy and the criteria for radiation in general relativity.
{\em Phys. Rev.} (2) 118 (1960) 1100-1104.

\bibitem{AM} L. J. Al\'{\i}as and J. M. Malacarne,
Spacelike hypersurfaces with constant higher order mean curvature in Minkowski space-time,
{\em J. Geom. Phys.} 41 (2002),  359-375.

\bibitem{ARS} L. J. Al\'{\i}as, A. Romero and M. S\'anchez, Spacelike hypersurfaces of constant mean curvature in certain spacetimes, {\em Nonlin. Analysis} 30 (1997), 655-661.

\bibitem{Ba1} R. Bartnik, The mass of an asymptotically flat manifold, {\em Comm. Pure Appl. Math.} 39 (1986), 5, 661-693.

\bibitem{Ba2} R. Bartnik,
Quasi-spherical metrics and prescribed scalar curvature,
{\em J. Differential Geom.} 37 (1993), 31-71.

\bibitem{Ba3} R. Bartnik,  Mass and 3-metrics of non-negative scalar curvature, {\em Proceedings of the International Congress of Mathematicians, Vol. II}, 231-240, Higher Ed. Press, Beijing, 2002.

\bibitem{Br} H. L. Bray,  Proof of the Riemannian Penrose inequality using the positive mass theorem. {\em J. Differential Geom.} 59 (2001), no. 2, 177-267.

\bibitem{BI} H. L. Bray and  K. Iga,
Superharmonic functions in $\mathbb R^n$ and the Penrose inequality in general relativity,
{\em Comm. Anal. Geom.} 10 (2002), 5, 999-1016.

\bibitem{BK} H. L. Bray and M. A. Khuri,  P.D.E.'s which imply the Penrose conjecture,  {\em Asian J. Math.} 15 (2011), no. 4, 557-610.

\bibitem{BL} H. L. Bray and D. L.  Lee,  On the Riemannian Penrose inequality in dimensions less than eight. {\em Duke Math. J.} 148 (2009), no. 1, 81-106.

\bibitem{Bre}  S. Brendle,  Rigidity phenomena involving scalar curvature, {\em arXiv:1008.3097}



\bibitem{BM} S. Brendle and F. Marques,  Recent progress on the Yamabe prolem, {\em arXiv:1010.4960}.

\bibitem{BMN} S. Brendle, F. Marques and A. Neves,  Deformations of the hemisphere that increase scalar
curvature, {\em Invent. Math.}, 2011, 185, 1, 175-197.


\bibitem{Ch}  P. Chru\'sciel, Boundary conditions at spatial infinity from a Hamiltonian point of view, {\em Topological properties and global structure of space-time}, 49–59, NATO Adv. Sci. Inst. Ser. B Phys., 138, Plenum, New York, 1986.

\bibitem{CH} P. Chru\'sciel and M. Herzlich,
The mass of asymptotically hyperbolic Riemannian manifolds,
{\em Pacific J. Math.} 212 (2003), no. 2, 231-264.

\bibitem{CN} P. Chru\'sciel and G. Nagy,  The mass of spacelike hypersurfaces in asymptotically anti-de Sitter space-times, {\em Adv. Theor. Math. Phys.} 5 (2001), 4, 697-754.

\bibitem{C} J. Corvino,  Scalar curvature deformation and a gluing construction for the Einstein constraint equations. {\em Comm. Math. Phys.} 214 (2000),  137-189.

\bibitem{dL} H. F. de Lima,   Spacelike hypersurfaces with constant higher order mean curvature in de Sitter space, {\em J. Geom. Phys.} 57 (2007), 967-975.

\bibitem{dLG} L. L. de Lima, and F. Gir\~ao,  Positive mass and Penrose type inequalities for asymptotically hyperbolic hypersurfaces, {\em  arXiv:1201.4991}.

\bibitem{FS} A. Freire and F. Schwartz,  Mass-capacity inequalities for conformally flat manifolds with boundary, {\em  arXiv:1107.1407}.

\bibitem{H1} M. Herzlich, Minimal surfaces, the Dirac operator and the Penrose inequality, {\em  S\'eminaire de Théorie Spectrale et G\'eom\'etrie}, Vol. 20, 2001-2002, 9-16.

\bibitem{H2} M. Herzlich,  Mass formulae for asymptotically hyperbolic manifolds, {\em AdS/CFT correspondence: Einstein metrics and their conformal boundaries}, 103-121, IRMA Lect. Math. Theor. Phys., 8, Eur. Math. Soc.,
    Z\"urich, 2005.

\bibitem{HW} L.-H Huang and D. Wu,  Hypersurfaces with nonnegative scalar curvature, {\em arXiv:1102.5749v2}.

\bibitem{HI} G. Huisken and T. Ilmanen,  The inverse mean curvature flow and the Riemannian Penrose inequality, {\em J. Differential Geom.} 59 (2001), no. 3, 353-437.

\bibitem{J} J. L. Jauregui,  Penrose-type inequalities with a Euclidean background, {\em  arXiv:1108.4042}.

\bibitem{KMS} M. A. Khuri, F. Marques and R. M. Schoen,
A compactness theorem for the Yamabe problem,
{\em J. Differential Geom.} 81 (2009), 1, 143-196.

\bibitem{L1} M.-K. G. Lam, The Graphs Cases of the Riemannian Positive Mass and Penrose Inequalities in All Dimensions, {\em arXiv:1010.4256}.

\bibitem{L2} M.-K. G. Lam,  The Graphs Cases of the Riemannian Positive Mass
and Penrose Inequalities in All Dimensions, {\em Duke thesis}.

\bibitem{LP} J. M. Lee and T. H.  Parker,  The Yamabe problem, {\em Bull. Amer. Math. Soc. (N.S.)} 17 (1987), 1, 37-91.

\bibitem{LZ} Y. Y. Liand L. Zhang,  Compactness of solutions to the Yamabe problem. II, {\em  Calc. Var. Partial
Differential Equations} 24 (2005),  2, 185-237.

\bibitem{M} B.  Michel,  Geometric invariance of mass-like asymptotic invariants, {\em J. Math. Phys.} 52 (2011), no. 5, 052504, 14 pp.



\bibitem{Ro} H. Rosenberg,  Hypersurfaces of constant curvature in space forms, {\em Bull. Sci. Math.} 117 (1993), 2, 211-239.

\bibitem{R} R. Reilly,  On the Hessian of a function and the curvatures of its graph, {\em  Michigan Math. J.}  20  (1973), 373-383.


\bibitem{S} R. Schoen,  Conformal deformation of a Riemannian metric to constant scalar curvature, {J. Differential Geom.} 20 (1984),  2, 479-495.


\bibitem{SY} R. Schoen and S.-T. Yau,  On the proof of the positive mass conjecture in general relativity. {\em Comm. Math. Phys.} 65 (1979), no. 1, 45-76.

\bibitem{Sc} F. Schwartz,  A volumetric Penrose inequality for conformally flat manifolds, {\em Ann. Henri
Poincar\'e}, 12 (2011), 67-76.

\bibitem{Sz} L. Szabados, Quasi-local energy-momentum and angular momentum in General Relativity,
{\em Living Rev.} 4 (2004), available at {\tt http://relativity.livingreviews.org/}

\bibitem{W} E. Witten,
A new proof of the positive energy theorem.
{\em Comm. Math. Phys.} 80 (1981), no. 3, 381-402.
\end{thebibliography}

\end{document}